\documentclass{jT}

\usepackage{amssymb,amsmath,latexsym,amsthm}
\usepackage[english]{babel}
\usepackage{hyperref}
\usepackage{ mathrsfs }
\usepackage{enumerate}
\usepackage{graphicx}

\usepackage{color}
\definecolor{darkblue}{rgb}{0,0,0.8}
\hypersetup{colorlinks,
            linkcolor=darkblue,
            anchorcolor=darkblue,
            citecolor=darkblue}
\let\orgautoref\autoref
\renewcommand{\autoref}[1]
        {\def\equationautorefname~##1\null{~(##1)\null}%
         \orgautoref{#1}}

\hypersetup{colorlinks=true}

\newtheorem{thm}{Theorem}[section]
\newtheorem{theorem}[thm]{Theorem}

\theoremstyle{definition}

\numberwithin{equation}{section}

\begin{document}

\title[The condition for an algebraic integer to be a Salem number]{The necessary and sufficient condition for an algebraic integer to be a Salem number}


\author{Dragan Stankov}
\address{Dragan Stankov\\
Katedra Matematike RGF-a Universiteta u Beogradu\\
11000 Beograd, Dju\v{s}ina 7\\
Serbia}
\email{dstankov@rgf.bg.ac.rs}





\maketitle


\begin{abstr}
We present a necessary and sufficient condition for a root greater than unity of a monic reciprocal polynomial of an even degree at least four,
with integer coefficients, to be a Salem number. We determine the probability of fulfillment the condition for an arbitrary power of the root.
\end{abstr}

\bigskip
\section{Introduction}
\label{intro}

A Salem number is a real algebraic integer $\tau > 1$ of degree at least four,
conjugate to $\tau ^{-1}$, all of whose conjugates, excluding $\tau$ and $\tau ^{-1}$, are unimodal i.e., lie on $|z| = 1$. 
The corresponding minimal polynomial $P(x)$ of degree $d$
of these numbers, called a Salem polynomial, is (self-)reciprocal, that is $x^d P(1/x) = P(x)$. Since $P(x)$ is self-reciprocal and irreducible it must have even degree.
It is well known \cite{Smy} that $\tau^n$ should also be a Salem number of degree $d$ for any natural $n$. Fractional parts of $\tau^n$ are dense in the unit interval $[0, 1]$, but are not uniformly distributed \cite{Ber,Sta}.
Salem numbers have appeared in quite different areas of mathematics (number theory, harmonic analysis, knot theory, etc.). Throughout, when we speak about a conjugate, the minimal
polynomial or the degree of an algebraic number we mean over the field of the rationals $\mathbb{Q}$.

In \cite{Vie} Vieira, extending a result of Lakatos and Losonczi \cite{LL}, presented a sufficient condition for a self-reciprocal
polynomial to have a fixed number of roots on the complex unit circle $U=\{z\in \mathbb{C}:|z|=1\}$.
Let $p(z) = a_d z^d+a_{d-1}z^{d-1}+\cdots+a_1z+a_0$ be a $d$-th degree self-reciprocal
polynomial. If the inequality
\begin{equation}\label{intro:IneqV}
|a_{d-l}|>\frac{1}{2}\left(\frac{d}{d-2l}\right)\sum_{k=0,k\ne l,d-l}^{d}|a_{k}|,\;\;l<d/2
\end{equation}
holds, then $p(z)$ has exactly $d - 2l$ roots on $U$ and these roots
are simple.
Here we present, in a sense, a result which lies in the opposite direction of a special case of this theorem.
Namely, we shall prove the following

\begin{theorem}\label{cha:Thm1}
A real algebraic integer $\tau > 1$ is a Salem number if and only if its minimal polynomial $P(x)$
is reciprocal of even degree $d\ge 4$, and there is $n\in \mathbb{N}$, $n\ge 2$ such that $\tau^n$ has the minimal polynomial $P_n(x)=1+a_{1,n}x+a_{2,n}x^2+\cdots+a_{d-1,n}x^{d-1}+x^d$, which is also reciprocal of degree $d$, and satisfies the condition
\begin{equation}\label{intro:Ineq1}
|a_{d-1,n}|>\frac{1}{2}\left(\frac{d}{d-2}\right)\left(2+\sum_{k=2}^{d-2}|a_{k,n}|\right).
\end{equation}
\end{theorem}
\noindent Notice that the condition \autoref{intro:Ineq1} is the special case when $l=1$ of the condition \autoref{intro:IneqV} applied to $P_n(x)$.

We present a method, easy for implementation, for the calculation of the coefficients of $P_n(x)$ starting with $P(x)$ without determination of its roots. We can use the companion matrix $C$ of a monic polynomial $P(x)=x^d+a_{d-1}x^{d-1}+a_{d-2}x^{d-2}+\cdots+a_0$ defined as
\[C=\left[ \begin{array}{cccccc}
0 &  1  & 0 & \cdots & 0 & 0\\
0 &  0  & 1 & \cdots & 0 & 0\\
\vdots &  \vdots  & \vdots & \vdots & \vdots & \vdots\\
0 &  0  & 0 & \cdots & 0 & 1\\
-a_0 & -a_1  & -a_2 & \cdots & -a_{d-2} & -a_{d-1}
\end{array}
\right]_{d\times d}
\]
It is well known \cite{LT}, \cite{LN} that $P(x)$ is the characteristic polynomial of $C$ so the root $\lambda$ of $P(x)$ is an eigenvalue of $C$. If $v$ is an
eigenvector of $C$ associated with $\lambda$ then $C^nv=C^{n-1}Cv=C^{n-1}\lambda v=\cdots=\lambda^n v$. Thus $C^n$ should have an eigenvalue $\lambda^n$ and the characteristic polynomial of $C^n$ must be $P_n(x)$, i.e. $P_n(x)=det(x I-C^n)$. It is easy to show that $v=[1\;\lambda\;\lambda^2\;\ldots\;\lambda^{d-1}]^T$.

Using this method we are able, for a Salem number $\tau$, to find at least one $n$ such that the minimal polynomial $P_n(x)$ of $\tau^n$ satisfies condition \autoref{intro:Ineq1}. In Table \autoref{tab:Pn} we present examples of Salem numbers and $n$ which we have found. The last example in the table is the root of Lehmer polynomial which is the smallest known Salem number.
We can notice that $n$ becomes large as $d$ increases. It would be interesting to find $n$ for all small Salem numbers in the Mossinghoff's list \cite{Mos}.

As shown in Table \autoref{tab:Pn} the relative frequency of $n$ such that the minimal polynomial $P_n(x)$ of $\tau^n$ satisfies \autoref{intro:Ineq1} significantly decreases when $d$ increases. One might ask what is the probability of fulfillment the condition \autoref{intro:Ineq1} for an arbitrary power of the root. We determined the exact value of the probability for $d=4, 6$ and we approximated the probability for $d=8, 10$.

\begin{theorem}\label{cha:ThmProb}
Let $\tau$ be a Salem number of degree $d$, $n\in \mathbb{N}$ and let $P_n(x)$ be the minimal polynomial of $\tau^n$. Let $p_d$ denotes the probability that coefficients of $P_n(x)$ satisfy \autoref{intro:Ineq1} when $n$ is randomly chosen. Then:

(a) $p_4$ is equal to $1/3$ and,

(b)
\begin{equation}\label{intro:Integral}
p_6=\frac{4}{\pi^2}\left[\int_{\arccos\frac{\sqrt{30}}{6}}^{\arccos\frac{\sqrt{19}-1}{6}}\left(\arccos\frac{-5-6\cos t}{6+6\cos t}-(\pi-t)\right)\mathrm{d}t+\right.\\
\end{equation}
\[\left.+\int_{\arccos\frac{\sqrt{19}-1}{6}}^{\arccos\frac{\sqrt{6}}{6}}\left(\arccos\frac{1-6\cos t}{6-6\cos t}-(\pi-t)\right)\mathrm{d}t\right]=0.0717258\ldots.\]


\end{theorem}

Furthermore, we have approximated the probabilities for $d=8$ and $d=10$ using a numerical method and have got $p_8\approx 0.012173$, $p_{10}\approx 0.0018$. These results suggest that $p_d$ decreases approximately five times when $d$ is increased by two.

If we observe coefficients of $P_n(x)$ as $n$ increases, we can notice some regularities which enable us to recognize the minimal polynomial of a Salem number. We present these regularities in the Theorem \autoref{cha:Thm2}.

\begin{table}
\caption{Salem number $\tau$ and first few $n$ such that the minimal polynomial $P_n(x)$ of $\tau^n$ satisfies (2)} 
\label{tab:Pn}       
\begin{tabular}{lrlll}
\hline\noalign{\smallskip}
 & d & $\tau$ & Coefficients & $n$: $P_n(x)$ satisfies \autoref{intro:Ineq1}  \\
\noalign{\smallskip}\hline\noalign{\smallskip}
1. & 4 &  1.72208381 & 1 -1 -1 & 9,13,16,17,20,24,27,31,35,38,42,45\\
2. & 6 &  1.50613568 & 1 -1 0 -1 & 14,16,35,37,54,65,67,86,116,144,157\\
3. & 8 &  1.28063816 & 1  0  0 -1 -1  & 72, 127, 163, 176\\
4. & 10 & 1.21639166 & 1  0  0  0 -1 -1 & 53\\
5. & 10 & 1.23039143 & 1  0  0 -1  0 -1 & 240\\
6. & 10 & 1.26123096 & 1  0 -1  0  0 -1 & 43, 80\\
7. & 10 & 1.17628082 & 1  1  0 -1 -1 -1 & 605\\
\noalign{\smallskip}\hline
\end{tabular}
\end{table}

\begin{theorem}\label{cha:Thm2}
Let $\tau > 1$ be a Salem number and let $P_n(x)=1+a_{1,n}x+a_{2,n}x^2+\cdots+a_{d-1,n}x^{d-1}+x^d$ be the minimal polynomial of $\tau^n$ for $n \in \mathbb{N}$. Then

(a) $\lim_{n\rightarrow \infty}\frac{a_{d-1,n+1}}{a_{d-1,n}}=\tau$, $\lim_{n\rightarrow \infty}\sqrt[n]{|a_{d-1,n}|}=\tau$,

(b) $\tau^n+\tau^{-n}\le |a_{d-1,n}|+d-2$,

(c) for $k=1,2,\ldots,d-3$
\begin{equation}\label{intro:Ineq3}
|a_{d-k-1,n}|\le \binom{d-2}{k}(|a_{d-1,n}|+d-2)+\binom{d-2}{k-1}+\binom{d-2}{k+1}.
\end{equation}
\end{theorem}

So if any of the conditions in Theorem \autoref{cha:Thm2} is not satisfied we can be sure that a root of $P(x)$ is not a Salem number. Theorem \autoref{cha:Thm2} explains the observation that the coefficients $a_{k,n}$ for $k=1,2,\ldots,d-1$ of $P_n(x)$ are approximately of the same magnitude, and that the central coefficient is usually slightly greater in modulus than a peripheral one. Examples for this are: $P_{100}$, 
$P_{200}$, showed in Table \autoref{tab:2}. The algorithm for calculating a root of $P(x)$ presented in (a) is known as Graeffe's method \cite{MP}.

If $P(x)$ is monic, reciprocal, with integer coefficients then $P_n(x)$ is a periodic sequence of polynomials if and only if $P(x)$ is the product of cyclotomic polynomials. In fact, if $P_n(x)$ is a periodic sequence, among these polynomials there are only finitely many distinct ones. Then the set of roots of these polynomials is also finite, and all the powers $\alpha,\alpha^2,\alpha^3,\ldots$ of a root $\alpha$ of $P(x)$ are in this set. Therefore for some $p$, $q$, $\alpha^p=\alpha^q$, $p\ne q$. Since $\alpha\ne 0$ it follows that $\alpha^{p-q}=1$. Vice versa, if $P(x)$ is the product of cyclotomic polynomials then all its roots are roots of 1 so the set of its powers is finite and the set of coefficients $a_{k,n}$ for $k=1,2,\ldots,d-1$, $n=1,2,\ldots$ of $P_n(x)$ is also finite. Thus $P_n(x)$ is a periodic sequence of polynomials.

\section{Proofs of Theorems}
\label{Proofs}

In order to prove Theorem \autoref{cha:Thm1} we shall use a theorem of Kronecker \cite[Theorem. 4.6.4.]{Ber}, which is a consequence of Weyl's theorems \cite{Bug}.
Suppose $\alpha = (\alpha_k)_{1\le k\le p}\in \mathbb{R}^p$ has the property that the real numbers
$1, \alpha_1, \ldots ,\alpha_p$ are $\mathbb{Q}$-linearly independent, and let $\mu$ denote an arbitrary vector
in $\mathbb{R}^p$, N an integer and $\varepsilon$ a positive real number. Then Kronecker's theorem states that there exists an integer
$n > N$ such that $\|n\alpha_k - \mu_k\| < \varepsilon$, $(k = 1,\ldots ,p)$ where
$\|x\| = \min \{ |x - m| : m \in \mathbb{Z}\}$ is the distance from $x$ to the nearest integer.

\begin{figure}[t]
  \includegraphics[width=0.7\textwidth]{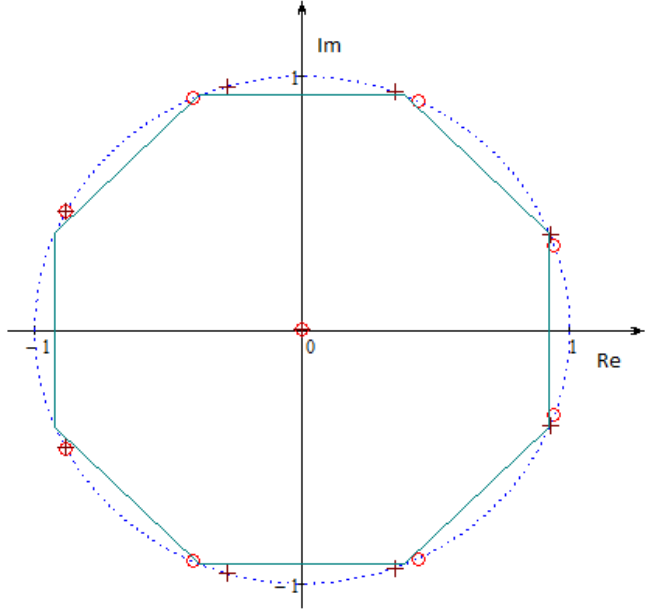}
\caption{If $\tau$ is the sixth Salem number in Table 1 of degree 10 then conjugates of $\tau^{43}$ (represented with $\circ$) and of $\tau^{80}$ (represented with $+$), whose minimal polynomials $P_{43}(x)$, $P_{80}(x)$ satisfy (2), are close to roots of $x^8+1$, vertices of the regular octagon.}
\label{fig:1}       
\end{figure}
%

\noindent
\textit{Proof of Theorem \autoref{cha:Thm1} Necessity}.
Suppose that $\tau > 1$ is a Salem number. The essence of the proof is to show that there is $n$ such that each of $d-2$ unimodal roots of $P_n(x)$ could be arbitrarily close to exactly one root of $x^{d-2}+1$ (see \cite[Lemma 2]{Zai}) and to show that then the coefficients of $P_n(x)$ will satisfy the condition \autoref{intro:Ineq1}. It is obvious that roots of $x^{d-2}+1$ are $\exp(\pm\frac{\pi+2j\pi}{d-2}i)$, $j=0,1,\ldots,d/2-2$.
We denote conjugates of $\tau$ by
\begin{equation}\label{intro:Conj}
\tau^{-1},\;\; \exp(\pm2i\pi\omega_1),\ldots ,\exp(\pm2i\pi\omega_{d/2-1}).
\end{equation}
Numbers $1$, $\omega_1,\ldots ,\omega_{d/2-1}$ are $\mathbb{Q}$-linearly independent \cite[Theorem 5.3.2.]{Ber}. According to the Kronecker's theorem consider
$(w_j)_{2\le j\le d/2}\in \mathbb{R}^{d/2-1}$ with $\mu = (\frac{1/2+0}{d-2}, \frac{1/2+1}{d-2}, \ldots ,\frac{1/2+d/2-2}{d-2})$. It is
clear that for every $\varepsilon >0$ there exists an arbitrarily large integer $n$ such that
\begin{equation}\label{intro:Kron}
|n\omega_{j}-\frac{1/2+j-1}{d-2}|<\varepsilon\;\;  (\mod 1)\;\; (j=1,2,\ldots,d/2-1).
\end{equation}
Since a coefficient of a polynomial is a continuous function of its roots, for every $\epsilon >0$ there exists an arbitrarily large integer $n$ such that the minimal polynomial
\begin{equation}\label{intro:epsi}
P_n(x)=(x-\tau^n)(x-\tau^{-n})(x^{d-2}+1+\sum_{j=1}^{d/2-2}\epsilon_{j}(x^{d-2-j}+x^{j})+\epsilon_{d/2-1}x^{d/2-1}),
\end{equation}
of the Salem number $\tau^n$ satisfies $|\epsilon_k|<\epsilon$, $k=1,\ldots,d/2-1$. We denote
\begin{equation}\label{intro:T}
-\tau^n-\tau^{-n}=T
\end{equation}
\newpage
\begin{eqnarray*}
P_n(x)&=&  (x^2+Tx+1)\cdot\\
&&\cdot(x^{d-2}+1+\sum_{j=1}^{d/2-2}\epsilon_{j}(x^{d-2-j}+x^{j})+\epsilon_{d/2-1}x^{d/2-1})\\
&=&x^{d}+1+(T+\epsilon_{1})(x^{d-1}+x)+(\epsilon_{2}+T\epsilon_{1}+1)(x^{d-2}+x^2)+\\
&&+\sum_{j=3}^{d/2-1}(\epsilon_j+T\epsilon_{j-1}+\epsilon_{j-2})(x^{d-j}+x^j)+\\
&&+(2\epsilon_{d/2-2}+T\epsilon_{d/2-1})x^{d/2}.
\end{eqnarray*}


Now we consider the coefficients of $P_n(x)$ to show they satisfy the condition \autoref{intro:Ineq1}. It is obvious that
$|a_{d-1,n}|=|T+\epsilon_{1}|\ge |T|-|\epsilon|$. We need to estimate
\begin{align*}
&\frac{1}{2}\left(\frac{d}{d-2}\right)(2+\sum_{k=2}^{d-2}|a_{k,n}|)=\\
=&\frac{1}{2}\left(\frac{d}{d-2}\right)(2+2|\epsilon_{2}+T\epsilon_{1}+1|+2\sum_{j=3}^{d/2-1}|\epsilon_j+T\epsilon_{j-1}+\epsilon_{j-2}|+\\
&+|2\epsilon_{d/2-2}+T\epsilon_{d/2-1}|)\\
\le&\frac{1}{2}\left(\frac{d}{d-2}\right)(2+2\epsilon_{2}+2|T|\epsilon_{1}+2+2\sum_{j=3}^{d/2-1}(\epsilon_j+|T|\epsilon_{j-1}+\epsilon_{j-2})+\\
&+2\epsilon_{d/2-2}+|T|\epsilon_{d/2-1})\\
\le&\frac{1}{2}\left(\frac{d}{d-2}\right)(2+2\epsilon+2|T|\epsilon+2+2\sum_{j=3}^{d/2-1}(\epsilon+|T|\epsilon+\epsilon)+2\epsilon+|T|\epsilon)\\
=&\frac{1}{2}\left(\frac{d}{d-2}\right)(4+4\epsilon+4(d/2-3)\epsilon+3|T|\epsilon+2(d/2-3)|T|\epsilon)\\
=&\frac{1}{2}\left(\frac{d}{d-2}\right)(4+(2d-8)\epsilon+(d-3)|T|\epsilon).
\end{align*}





\noindent
So the condition \autoref{intro:Ineq1} will be satisfied if
\[|T|-|\epsilon|>\frac{1}{2}\left(\frac{d}{d-2}\right)(4+(2d-8)\epsilon+(d-3)|T|\epsilon),\]
\noindent
which is equivalent to
\begin{equation}\label{intro:Ineq2}
\frac{(2d-4)|T|-4d}{(d^2-3d)|T|+2d^2-6d-4}>\epsilon,\;\; d\ge 4.
\end{equation}
Since $|T|=\tau^n+\tau^{-n}$ tends to $\infty$ as $n\rightarrow\infty$ it is obvious that the left side of \autoref{intro:Ineq2} tends to $D:=\frac{2d-4}{d(d-3)}$
as $n\rightarrow\infty$. The determination of $n$ such that coefficients of $P_n(x)$ satisfies \autoref{intro:Ineq1} has to be done in following four steps:
\begin{enumerate}[i]
    \item we choose $\epsilon$ such that $D>\epsilon>0$;
    \item we choose an integer $N$ such that \autoref{intro:Ineq2} will be fulfilled for all $n\ge N$;
    \item we chose an $\varepsilon>0$ such that if each of $d-2$ unimodal roots of a $P_n(x)$ is at the distance $<\varepsilon$ in modulus of exactly one root of $x^{d-2}+1$ then $|\epsilon_k|<\epsilon$, $k=1,\ldots,d/2-1$ is fulfilled in \autoref{intro:epsi};
    \item we chose $n\ge N$ such that \autoref{intro:Kron} is fulfilled.
\end{enumerate}

\noindent
\textit{Sufficiency}. Suppose that $\tau > 1$ is a real algebraic integer with conjugates $\tau_1 = \tau$, $\tau_2,\ldots,\tau_d$ over $\mathbb{Q}$ such that $\tau^n$ has the minimal polynomial $P_n(x)$ which is also reciprocal of degree $d$, and satisfies the condition \autoref{intro:Ineq1}.
If $\tau$ is a conjugate of $\tau'$ then $\tau^n$ is a conjugate of $\tau'^n$. Since the minimal polynomial $P_n(x)$ of $\tau^n$ is of degree $d$ so $\tau_1^n$, $\tau_2^n,\ldots,\tau_d^n$ must be different numbers and their product has to be 1 because $P_n(x)$ is monic and reciprocal. The polynomial $P_n(x)$ satisfies the condition \autoref{intro:Ineq1} so it satisfies the condition \autoref{intro:IneqV} of Vieira's theorem where $l=1$. According to the theorem there are $d-2$ roots of $P_n(x)$ on the boundary of the unit disc $|z|=1$. Since they occur in conjugate complex pairs their product is equal to 1. It follows that $\tau^{-n}$ should be a conjugate of $\tau^n$ which allow us to conclude that $\tau^{n}$ is a Salem number. If $|\tau'^n|=1$ then $|\tau'|=1$ thus it follows that there are $d-2$ conjugates of $\tau$ on the boundary of the unit disc. Finally, in the same manner as for $\tau^n$, we conclude that $\tau$ is also a Salem number.\qed

\noindent
\textit{Proof of Theorem \autoref{cha:ThmProb} (a)}
If we use \autoref{intro:Conj} and denote $D:=\tau^n+1/\tau^n$ ($d=4$) we have
\[P_n(x)=(x^2-Dx+1)(x^2-2\cos(2\pi n\omega_1)x+1).\]
We denote $2\pi \{n \omega_1\}$ by $\theta_1$ and $2\cos(\theta_1)$ by $s_1$ where $\{\cdot\}$ denotes the fractional part. Since $n\omega_1$ is uniformly distributed modulo one $\theta_1$ is uniformly distributed on $[0,2\pi]$. For $d=4$ the condition \autoref{intro:Ineq1} is reduced to $|a_{3,n}|>2+|a_{2,n}|$. Since $P_n(x)=(x^2-Dx+1)(x^2-s_1x+1)$ the condition becomes
\begin{equation}\label{intro:Ineq1deg4}
|-D-s_1|>2+|Ds_1+2|.
\end{equation}
From the definition of $D$ it is obvious that $D\rightarrow \infty$ when $n\rightarrow \infty$. Since $|s_1|\le 2$ we have $D+s_1\rightarrow \infty$ so that $|-D-s_1|=|D+s_1|$ is equal, for every sufficiently large $n$, to $D+s_1$. Finally \autoref{intro:Ineq1deg4} becomes $D+s_1>2+|Ds_1+2|$ i.e. $D+s_1-2>Ds_1+2>-D-s_1+2$. Solving this double inequality for $s_1$ we get
\[\frac{4-D}{D-1}<s_1<\frac{D}{D+1}.\]
When $D$ tends to infinity we obtain $-1<s_1<1$ i.e. $-1/2<\cos \theta_1<1/2$. It follows that $\pi/3<\theta_1<2\pi/3$ or $4\pi/3<\theta_1<5\pi/3$ so that the probability has to be $p_4=\frac{2\pi/3}{2\pi}=\frac{1}{3}$.

\textit{(b)}
Using \autoref{intro:Conj} with $d=6$  and the definition of $D$ we have
\[P_n(x)=(x^2-Dx+1)(x^2-2\cos(2\pi n\omega_1)x+1)(x^2-2\cos(2\pi n\omega_2)x+1).\]
We denote $\theta_1:=2\pi \{n \omega_1\}$, $\theta_2:=2\pi \{n \omega_2\}$.
Coefficients of $P_n(x)$ depends only on real parts of unimodal roots so that we can chose the complex conjugates from the upper half (complex) plane. Thus we define
\begin{equation}\label{proof:ti}
t_i = \left\{ \begin{array}{ll}
         \theta_i & \mbox{if $\theta_i\in(0,\pi)$};\\
        2\pi-\theta_i & \mbox{if $\theta_i\in(\pi,2\pi)$}.\end{array} \right.,\;i=1,2.
\end{equation}

\begin{figure}[t]
  \includegraphics[width=0.7\textwidth]{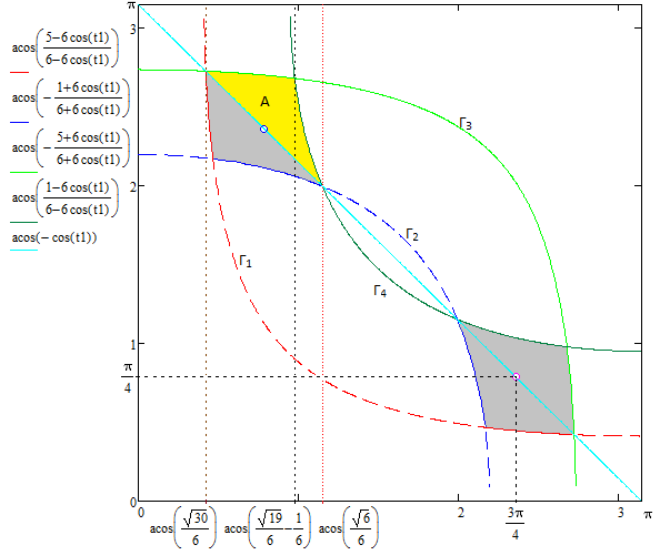}
\caption{The event T that a power of a Salem number of degree 6 has the minimal polynomial which satisfies the condition (1.2) is shaded in the figure. It consists of four congruent curve-triangles, each of them has the same area $A$ which is equal to the definite integral. Thus the probability of $T$ is $4A/{\pi}^2=0.0717258\ldots$. }
\label{fig:2}       
\end{figure}

Since $n\omega_1$, $n\omega_2$ are uniformly distributed modulo one $\theta_1$, $\theta_2$ are uniformly distributed on $[0,2\pi]$ and $t_1$, $t_2$ are uniformly distributed on $[0,\pi]$.
We denote
\begin{equation}\label{proof:st}
s_1:=2\cos(t_1), s_2:=2\cos(t_2).
\end{equation}
For $d=6$ the condition \autoref{intro:Ineq1} is reduced to $|a_{5,n}|>\frac{6}{8}(2+2|a_{4,n}|+|a_{3,n}|)$. Since
\[P_n(x)=(x^2-Dx+1)(x^2-s_1x+1)(x^2-s_2x+1)\] the condition becomes
\begin{equation}\label{intro:Ineq1deg6}
|-D-s_1-s_2|>\frac{6}{8}(2+2|Ds_1+Ds_2+s_1s_2+3|+|-2D-2s_1-2s_2-Ds_1s_2|.
\end{equation}
The main idea of the proof is to determine the region $S$ in $s_1Os_2$ plane such that every point $(s_1,s_2)\in S$ satisfies \autoref{intro:Ineq1deg6}. Since $D\rightarrow \infty$ when $n\rightarrow \infty$, $|s_1|\le 2$, $|s_2|\le 2$ we conclude that the left side in \autoref{intro:Ineq1deg6} $|-D-s_1-s_2|=|D+s_1+s_2|$ is equal, for every sufficiently large $n$, to $D+s_1+s_2$.
We can find the boundary of $S$ if we replace $>$ in \autoref{intro:Ineq1deg6} with $=$ and if we replace both $|\;|$ on the right side with $\pm(\;)$. There are four possibilities for replacing so we get four equations which we solve for $s_2$. We get rational functions $s_2=f_i(D,s_1)$ which tends to $s_2=F_i(s_1)$ when $D\rightarrow\infty$, $i=1,2,3,4$:
\[f_1(D,s_1)=\frac{10D+10s_1-6Ds_1-24}{6D+6s_1-3Ds_1-10},\;F_1(s_1)=\frac{10-6s_1}{6-3s_1},\]
\[f_2(D,s_1)=-\frac{2D+2s_1+6Ds_1+24}{6D+6s_1+3Ds_1+2},\;F_2(s_1)=-\frac{2+6s_1}{6+3s_1},\]
\[f_3(D,s_1)=-\frac{10D+10s_1+6Ds_1+12}{6D+6s_1+3Ds_1+10},\;F_3(s_1)=-\frac{10+6s_1}{6+3s_1},\]
\[f_4(D,s_1)=\frac{2D+2s_1-6Ds_1-12}{6D+6s_1-3Ds_1-2},\;F_4(s_1)=\frac{2-6s_1}{6-3s_1}.\]
The boundary of $S$ consists of parts of graphs of $F_i(s_1)$. We have to find intersection points of these graphs. Therefore we solve four equations:
\[F_1(s_1)=F_2(s_1)\Rightarrow s_1=1/3\pm\sqrt{19}/3,\]
\begin{equation}\label{proof:intsec}
F_1(s_1)=F_3(s_1)\Rightarrow s_1=\pm\sqrt{30}/3,
\end{equation}
\[F_4(s_1)=F_2(s_1)\Rightarrow s_1=\pm\sqrt{6}/3,\]
\[F_4(s_1)=F_3(s_1)\Rightarrow s_1=-1/3\pm\sqrt{19}/3.\]
We have to determine the area of the region $T$ in $t_1Ot_2$ plane such that for every point $(s_1,s_2)\in S$ there is unique $(t_1,t_2)\in T$ where
\begin{equation}\label{proof:ts}
t_1=\arccos(s_1/2),\; t_2=\arccos(s_2/2),
\end{equation}
using \autoref{proof:st}. The ratio of the area of $T$ to the area of all possible values $(t_1,t_2)$, i.e. $\pi^2$, is equal to the probability $p_6$. Since $s_2=F_i(s_1)$ it follows that $t_2=\arccos(F_i(2\cos(t_1))/2)=:G_i(t_1)$ using \autoref{proof:st}. For the determination of the area of $T$ it is convenient to show that $T$ has reflection symmetry across the line $t_2=\pi-t_1$.
Let the graph of $t_2=G_i(t_1)$ be $\Gamma_i$. We claim that $\Gamma_1$ can be obtained by reflecting of $\Gamma_3$ about the line $t_2 = \pi-t_1$ i.e. if $(t_1,t_2)\in\Gamma_1$ then $(\pi-t_2,\pi-t_1)\in\Gamma_3$ (see Figure \autoref{fig:2}). Indeed, if $t_2=G_1(t_1)$ then
\begin{eqnarray*}
G_3(\pi-t_2)&=&\arccos(F_3(2\cos(\pi-t_2))/2)\\
&=&\arccos(F_3(-2\cos(t_2))/2)\\
&=&\arccos(F_3(-2\cos(G_1(t_1)))/2)\\
&=&\arccos(F_3(-2F_1(2\cos(t_1))/2)/2)\\
&=&\arccos(F_3(-F_1(s_1))/2)\\
&=&\arccos(-s_1/2)\\
&=&\arccos(-\cos(t_1))\\
&=&\arccos(\cos(\pi-t_1))\\
&=&\pi-t_1.
\end{eqnarray*}

In the same manner we can show that $\Gamma_2$ is a reflection of $\Gamma_4$ in the line $t_2 = \pi-t_1$. Therefore $T$ consists of four congruent curve-triangles, each of them has the same area $A$ (see Figure \autoref{fig:2}). If we bring to mind the intersection points \autoref{proof:intsec} and formulas \autoref{proof:ts} we find out the intersection points of graphs $\Gamma_i$, $i=1,2,3,4$ which are the limits of two definite integrals which occur in \autoref{intro:Integral}. We conclude that $A$ is equal to sum of these integrals (see Figure \autoref{fig:2}) and that $p_6=4A/\pi^2$ as it is claimed.
\qed

If we use the same method for the determination of $p_8$, $p_{10}$ etc. it requires multiple definite integrals applied on the regions with complicated boundaries. Thus it is much convenient to use a numerical approach. For each pair of conjugate complex roots of a Salem polynomial we define a variable $t_i\in(0,\pi)$, as in \autoref{proof:ti} and $s_i$ as in \autoref{proof:st} $i=1,2,\ldots,H$ where we denoted $(d-2)/2$ by $H$. Let $m\in \mathbb{N}$ and let $0=t_{i,0},t_{i,1},\ldots,t_{i,m}=\pi$, $i=1,2,\ldots,H$ be nodes arranged consecutively with equal spacing $h=\pi/m$. Starting from
\[P_n(x)=(x^2-Dx+1)\prod_{i=1}^{H}(x^2-2\cos(t_i)x+1)\]
we calculate the coefficients of $P_n(x)$ which obviously depend on $D$, $t_i$ so that there are the functions $A_{k,n}$ such that
\[a_{k,n}=A_{k,n}(D;t_1,t_2,\ldots,t_{H}),\; k=1,2,\ldots,d-1.\]
For $D$ fixed and for each $H$-tuple $(t_{1,j_1},t_{2,j_2},\ldots,t_{H,j_H})$ we calculate
\[a_{k,n}=A_{k,n}(D;t_{1,j_1},t_{2,j_2},\ldots,t_{H,j_H}),\; j_i=0,1,\ldots,m,\]
and replace them into the condition \autoref{intro:Ineq1}. The number $N_c$ of all $H$-tuples, i.e. of all points of $\pi^H$, which satisfy this condition, divided with $(m+1)^H$, the number of all $H$-tuples, approximates $p_d$. If we take a large $D=10^9$ and a small $h\ge 0.002$ we get $p_8\approx 0.012173$, $p_{10}\approx 0.0018$. Since there are four nested loops the calculation of $p_{10}$ requires much CPU time. Thus it was necessary to improve our programm. We use the fact that all $H$-tuples which satisfy \autoref{intro:Ineq1} are close to the point $P(\pi/H,3\pi/H,\ldots,(d-3)\pi/H)$ or to $H!$ points obtained by permuting the coordinates of $P$, because these coordinates are the arguments of the roots of $x^{d-2}+1$. Therefore to get $N_c$ we have to check and count only points in a small region around the $P$ and then to multiply the number of them by $H!$. Executing the program with a small $D$ we have got less probability than with a large one. It suggest us that the convergence of $p_d^{n_0}$ to $p_d$ is from below.

We have also verified $p_4$ and $p_6$ statistically. For the first Salem number in the Table \autoref{tab:Pn} of degree 4 we have found that if $1\le n\le 300$ then the coefficients of $P_n(x)$ satisfy \autoref{intro:Ineq1} $98$ times: for
9, 13, 16, 17, 20, 24, 27, 31, 35, 38, 42, 45, 46, 49, 53, 56, 57, 60, 64, 67, 68, 71, 75, 78, 79, 82, 86, 89, 93, 97, 100,
104, 107, 108, 111, 115, 118, 122, 126, 129, 130, 133, 137, 140, 141, 144, 148, 151, 155, 159,162, 166, 169, 170, 173, 177, 180, 181, 184, 188, 191, 192, 195, 199,
202, 203, 206, 210, 213, 217, 221, 224, 228, 231, 232, 235, 239, 242, 243, 244, 246, 250, 253, 254, 257, 261, 264, 265, 268, 272, 275, 279, 283, 286, 290, 293, 294, 297,
so that the relative frequency is $98/300\approx 0.33$.

For the second Salem number in the Table \autoref{tab:Pn} of degree 6 we have found that if $101\le n\le 300$ then the event that $P_n(x)$ satisfies \autoref{intro:Ineq1} occurs fourteen times: for $n=116$, $144$, $157$, $167$, $187$, $195$, $206$, $225$, $238$, $246$, $257$, $276$, $287$, $295$ so that the relative frequency is $14/200=0.07$.
If $1001\le n\le 1200$ then $P_n(x)$ satisfy \autoref{intro:Ineq1} sixteen times: for $n=1001$, $1029$, $1031$, $1039$, $1050$, $1052$, $1063$, $1080$, $1082$, $1101$, $1103$, $1120$, $1131$, $1133$, $1152$, $1182$ with the relative frequency $16/200=0.08$.

\begin{table}
\caption{Coeficients of $P_{43}(x)$, $P_{80}(x)$ which satisfy (2) and of $P_{100}(x)$, $P_{200}(x)$ which do not, where $P(x)$ is the minimal polynomial of the sixth Salem number in Table 1}
\label{tab:2}       
\begin{tabular}{rrrrr}
\hline\noalign{\smallskip}
$P_{43}(x)$ & $P_{80}(x)$ & $P_{100}(x)$ & $P_{200}(x)$ \\ 
\noalign{\smallskip}\hline\noalign{\smallskip}
1 & 1 & 1 & 1 \\
-21586 & -115763027 & -12007769482 & -144186527874521531930 \\
3611 & 23986075 & 29164508197 & 415053787386817223949 \\
688 & -39926871 & -18134706516 & -542626204385602820124 \\
5418 & 20167702 & -25180138718 & 625113687841885675082 \\
-6193 & 4830711 & 52256753515 & -707660656174865919717 \\
\noalign{\smallskip}\hline
\end{tabular}
\end{table}

\noindent
\textit{Proof of Theorem \autoref{cha:Thm2}}
(a) Since $\tau^n$ is a Salem number $P_n(x)$ has to be monic, reciprocal polynomial of even degree. Using the notation \autoref{intro:Conj} for conjugates of $\tau$ and Vieta's formulae we have
\[\frac{a_{d-1,n+1}}{a_{d-1,n}}=\frac{\tau^{n+1}+\tau^{-n-1}+\sum_{k=1}^{d/2-1}\left(e^{2(n+1)i\pi\omega_k}+e^{-2(n+1)i\pi\omega_k}\right)}
{\tau^{n}+\tau^{-n}+\sum_{k=1}^{d/2-1}\left(e^{2ni\pi\omega_k}+e^{-2ni\pi\omega_k}\right)}.\]
If we divide the enumerator and the denominator with $\tau^{n}$ then it is obvious that we obtain the enumerator which tends to $\tau$ and the denominator which tends to $1$, as $n\rightarrow\infty$. If we write
\[|a_{d-1,n}|=\tau^{n}|1+\tau^{-2n}+\tau^{-n}\sum_{k=1}^{d/2-1}\left(e^{2ni\pi\omega_k}+e^{-2ni\pi\omega_k}\right)|,\]
then we can conclude immediately that $\lim_{n\rightarrow \infty}\sqrt[n]{|a_{d-1,n}|}=\tau$.

\noindent
(b) Since $a_{d-1,n}=-\tau^{n}-\tau^{-n}-\sum_{k=1}^{d/2-1}\left(e^{2ni\pi\omega_k}+e^{-2ni\pi\omega_k}\right)$
it is obvious that $\tau^n+\tau^{-n}\le |a_{d-1,n}|+d-2$.

\noindent
(c) From \autoref{intro:Conj} we have
\[P_n(x)=(x-\tau^n)(x-\tau^{-n})\prod_{k=1}^{d/2-1}\left[\left(x-e^{2ni\pi\omega_k}\right)\left(x-e^{-2ni\pi\omega_k}\right)\right],\]
\[P_n(x)=(x^2-(\tau^n+\tau^{-n})x+1)\prod_{m=1}^{d-2}\left(x-e^{(-1)^m2ni\pi\omega_{\left\lceil \frac{m}{2}\right\rceil}}\right)\]
where $k=\left\lceil \frac{m}{2}\right\rceil$. If we denote the product by $Q_n(x)=1+b_{1,n}x+b_{2,n}x^2+\cdots+b_{d-3,n}x^{d-3}+x^{d-2}$ then we can see that $b_{k,n}$, $k=1,2,\ldots,d-3$ is the sum of $\binom{d-2}{k}$ summands where each of them is of modulus $1$ so that $|b_{k,n}|\le \binom{d-2}{k}$. By expanding $P_n(x)=(x^2-(\tau^n+\tau^{-n})x+1)Q_n(x)$ we get for $k=1,2,\ldots,d-3$
\[a_{d-k-1,n}= -(\tau^n+\tau^{-n})b_{k,n}+b_{k-1,n}+b_{k+1,n}.\]
We can conclude that \autoref{intro:Ineq3} is valid using (b)
\[|a_{d-k-1,n}|\le \binom{d-2}{k}(|a_{d-1,n}|+d-2)+\binom{d-2}{k-1}+\binom{d-2}{k+1}.\]
\qed

\bigskip


\end{document}